\documentclass[11pt]{article}

\usepackage[T1]{fontenc}
\usepackage[utf8]{inputenc}
\usepackage{lmodern}
\usepackage{microtype}

\usepackage{amsmath,amssymb,amsfonts,mathtools}
\usepackage{dsfont}

\usepackage{graphicx}
\usepackage{tikz, pgfplots}
\usetikzlibrary{positioning}
\usetikzlibrary{3d}
\pgfplotsset{compat=1.18}

\usepackage{titlesec}

\titleformat{\section}
  {\normalfont\Large\bfseries}{\thesection.}{1em}{}

\titleformat{\subsection}
  {\normalfont\large\bfseries}{\thesubsection.}{1em}{}

\titleformat{\subsubsection}
  {\normalfont\normalsize\bfseries}{\thesubsubsection.}{1em}{}

\usepackage{subfig}

\usepackage{blindtext}
\usepackage{authblk} 

\usepackage[hidelinks]{hyperref}
\usepackage{orcidlink} 

\usepackage[backend=biber]{biblatex}
\newcommand{\orcidID}[1]{\orcidlink{#1}}

\newcommand{\keywords}[1]{%
  \begingroup
  \def\and{, }%
  \noindent\textbf{Keywords: }#1\par
  \endgroup}

\usepackage{amsthm}
\makeatletter
\@ifundefined{c@theorem}{%
  \theoremstyle{plain}
  \newtheorem{theorem}{Theorem}[section]
  \newtheorem{lemma}[theorem]{Lemma}

  \theoremstyle{definition}
  \newtheorem{definition}[theorem]{Definition}
  
  \theoremstyle{remark}
  
}{}
\makeatother

\def\R{{\mathbb R}}

\def\N{{\mathbb N}}

\def\P{\mathbb P} 

\def\X{\mathbb X}
\def\M{\mathbb M}
\def\L{\mathbb L}

\def\car#1{{\mathbf 1}}
\def\carp#1{{\mathbf 1}_{#1}} 

\newcommand{\esp}[1]{{\mathbb E}\left[#1\right]}      
\newcommand{\espi}[2]{{\mathbb E}_{#1}\left[{#2}\right]}    

\def\H{{\mathcal H}}

\def\B{{\mathcal B}} 
\def\S{\mathbb S} 
\def\cyl{\mathcal S}

\def\LL{\mathbf{L}}

\def\d{\operatorname{d}}  
\def\/{\,|\,}   
\newcommand{\Dom}{\operatorname{Dom}}

\def \DTV{\operatorname{d}_{\operatorname{TV}}} 
\def \DW{\operatorname{d}_{\operatorname{W}}} 
\def\Lip{\operatorname{Lip}}
\def\Leb{\ell} 

\def\NC{{\mathfrak N}} 
\def\FNC{{\mathfrak N \textsuperscript{f}}} 

\def\G{{\mathcal G}}

\def\pGlaub{\mathcal P}
\def\lGlaub{\mathcal L}

\newcommand{\step}[1]{\noindent{\small\textsc{Step }{#1}. }}

\def\eqdist{\overset{\mathcal{D}}{=}} 
\def\cvdist{\overset{\mathcal{D}}{\to}} 

\title{Quantitative Limit Theorems for Cox-Poisson and Cox-Binomial Point Processes}

\author[1]{Hamza Adrat\orcidID{0009-0004-8155-3260}}
\author[2]{Laurent Decreusefond\orcidID{0000-0002-8964-0957}}
\affil[1]{University Mohammed VI Polytechnic, Ben Guerir, Morocco\\
  \texttt{Hamza.ADRAT@um6p.ma}}
\affil[2]{Telecom Paris, Paris, France\\
  \texttt{Laurent.DECREUSEFOND@telecom-paris.fr}}

\date{} 

\begin{document}
\maketitle

\begin{abstract}
This paper establishes quantitative limit theorems for two classes of Cox point processes, quantifying their convergence to a Poisson point process (PPP). 
We employ Stein's method for PPP approximation, leveraging the generator approach and the Stein-Dirichlet representation formula associated with the Glauber dynamics. 
First, we investigate a Cox-Poisson process constructed by placing one-dimensional PPPs on the lines of a Poisson line process in $\R^2$. 
We derive an explicit bound on the convergence rate to a homogeneous PPP as the line intensity grows and the point intensity on each line diminishes. 
Second, we analyze a Cox-Binomial process on the unit sphere $\S^2$, modeling a system of satellites. This process is generated by placing PPPs on great-circle orbits, 
whose positions are determined by a Binomial point process. For this model, we establish a convergence rate of order $O(1/n)$ to a uniform PPP on the sphere, 
where $n$ is the number of orbits. The derived bounds provide precise control over the approximation error in both models, 
with applications in stochastic geometry and spatial statistics.
\end{abstract}


\pagenumbering{arabic}\setcounter{page}{1}

\section{Introduction}
\label{sec:introduction}

The Poisson point process (PPP) is a cornerstone of stochastic geometry and spatial statistics, serving as a fundamental model for phenomena characterized 
by complete spatial randomness. Its analytical tractability and elegant properties make it a natural benchmark for more complex spatial patterns. 
A common theme in the theory of point processes is the emergence of Poissonian behavior in the limit of superpositions or transformations of other, 
often non-Poisson, processes. While many results establish such convergence in law, a crucial aspect for practical applications is to quantify 
the rate of this convergence.

Stein's method, originally developed for normal approximation by Charles Stein \cite{Stein1972}, 
offers a powerful and versatile framework for obtaining explicit error bounds in probabilistic approximations. 
The generator approach to Stein's method is particularly well suited for point process settings, 
built upon the theory of spatial birth-and-death processes \cite{Preston1975}. This approach relies on characterizing the target distribution (in our case, the PPP) 
via a differential operator of a stationary Markov process, typically the Glauber dynamics. B
y solving the associated Stein equation through its semigroup, one can obtain integral representations for the difference between expectations, 
leading to explicit error bounds. This technique has been successfully adapted to various settings, 
including infinite-dimensional spaces and Brownian approximations \cite{Shih2011, CoutinDecreusefond2013}.

\subsection*{Our contribution}
This paper applies the Stein-generator framework to derive quantitative limit theorems for two important classes of Cox point processes, 
which are doubly stochastic processes where the intensity measure itself is random.
\begin{enumerate}
  \item The first model is a Cox-Poisson process built upon a Poisson line process in the Euclidean plane, where points are scattered along these random lines. 
  Such processes are instrumental in the modeling of wireless and vehicular networks \cite{DhillonChetlur2020}. 
  \item The second model is a Cox-Binomial process on the unit sphere, designed to represent satellite constellations where satellites are distributed 
  along great-circle orbits determined by a finite number of parent bodies.
\end{enumerate}

\section{Notations}
\label{sec:notations}
The following notions concern point process theory and come essentially from \cite{Daley2003, Decreusefond2022}.
We consider a locally compact metric space $\X$ endowed with its Borel space $\B(\X)$ and a measure $m$ on $\X$ (which may not be necessarily diffuse).
We denote $\NC_{\X}$ the space of configurations on $\X$ and $\FNC_{\X}$ designs the space of finite configurations on $\X$.
In particular, we denote by $\NC_{\X}^{(N)}$ the space of configurations on $\X$ with exactly $N$ points.
We endow $\NC_{\X}$ with $\B \left( \NC_{\X} \right)$ defined as the smallest $\sigma$-algebra on $\NC_{\X}$ such that $\phi \in \NC_{\X} \mapsto \phi(A)$
is measurable for any $A \in \B \left( \NC_{\X} \right)$. The restriction of $\B \left( \NC_{\X} \right)$ to $\FNC_{\X}$ is denoted $\B \left( \FNC_{\X} \right)$.
\noindent
For any $\phi, \omega$ in $\NC_{\X}$ and $x \in \X$, we will use the following notations:
\begin{equation*}
	\phi \oplus x =
	\left\{ \begin{array}{l}
		\phi, \text { if } x \in \phi, \\
		\phi \cup\{x\}, \text { if } x \notin \phi.
	\end{array} \right.
\end{equation*}
And,
\begin{equation*}
	\phi \ominus x=
	\left\{ \begin{array}{l}
		\phi \setminus \{x\}, \text { if } x \in \phi, \\
		\phi, \text { if } x \notin \phi.
	\end{array} \right.
\end{equation*}
Similarly, we will write $\phi \oplus \omega$ instead of $\phi \cup \omega$, $\phi \ominus \omega$ instead of $\phi \setminus \omega$, 
and $\phi(\X)$ will also be denoted as $|\phi|$.

\noindent
For $p \ge 1$, and a reference measure $\nu$ on $\left( \NC_{\X}, \B(\NC_{\X}) \right)$, we consider the space
\begin{equation*}
	\LL^p (\NC_{\X} \to \R, \, \nu) = \{ F \colon \NC_{\X} \to \R \; : \; \int_{\NC_{\X}} \left|F(\phi)\right|^p \nu(\d \phi) < \infty \}.
\end{equation*}

\noindent
We define a marked point process as a pair $\Psi = (\Phi, M)$, where $\Phi$ is a point process and each point $x \in \Phi$ is associated with a random variable $M(x)$, 
called a mark, taking values in a measurable space $\M$ (the marks space). 
The marks $M(x)$ can take various forms, such as integers, real numbers, geometric objects, or even other point processes.
In this case, $\Psi$ is a point process on $\X \times \M$.

\noindent
Finally, we denote $\cyl(\X)$ the space of functionals $F \colon \NC_{\X} \to \R$ such that there exists a compact subset $K$ of $\X$,
for any point process $\Phi$ on $\X$, we have $F(\Phi) = F\left(\Phi \cap K \right)$. Such a functional $F$ is called a cylindrical functional.

\section{Preliminaries}
\label{sec:preliminaries}

\subsection{Poisson and Binomial point processes}

\begin{definition}[Poisson point process]
  \label{def:poisson_process}
  Let $\sigma$ be a Radon measure on $\X$. The Poisson point process (PPP) $\Phi$ with intensity measure $\sigma$ is defined 
  as the unique point process on $\X$ with intensity measure $\sigma$ such that, for any disjoint relatively compact subsets 
  $K_1, K_2$, the random variables $\Phi\left(K_1\right)$ and $\Phi\left(K_2\right)$ are independent.

  \noindent
  Moreover, if there exists $\lambda \in[0,+\infty)$ such that $\sigma(\d x)=\lambda m(\d x)$, 
  then $\Phi$ is said to be homogeneous with intensity $\lambda$.
\end{definition}

\noindent
Let $\Phi$ be a marked point process, we denote by $\Phi^{0}$ the unmarked point process on $\X$, and we suppose that $\M$ is	equipped with measure $\nu$.
If $\Phi^{0}$ is a PPP with intensity $\sigma$, and the marks are stochastically independent, 
then the marked point process $\Phi$ is also a PPP on $\X \times \M$ with intensity $\sigma \otimes \nu$, called a marked Poisson point process (marked PPP).

\noindent
Similarly, we define a Binomial point process (BPP) and a marked BPP as follows.
\begin{definition}[Binomial point process]
	\label{def:binomial_process}
	Let $\mu$ be a probability measure on $\X$ and $N \in \N^*$. 
	A Binomial point process with parameter $N$ and supported by $\mu$ has exactly $N$ points drawn independently according to $\mu$.
\end{definition}
And if $\Phi^0$ is a BPP with parameter $N$ and supported by $\mu$, being the unmarked point process of $\Phi$, then $\Phi$ is also a BPP with 
parameter $N$ and supported by $\mu \otimes \nu$, where $\nu$ is the measure of the marks space $\M$.

\noindent
One of the fundamental formulas for PPPs is the following (see~\cite{Decreusefond2022}).
\begin{theorem}[Campbell-Mecke formula for a PPP]
  \label{th:mecke_formula_PPP}
	Let $\Phi$ be a PPP on $\X$ with intensity measure $\sigma$. 
	Then, for any measurable function $F \colon \X \times \NC_{\X} \to \R$ such that
	\begin{equation*}
		\esp{\int_{\X} F(x, \Phi) \sigma(\d x)} < \infty,
	\end{equation*}
	we have
	\begin{equation*}
		\esp{\sum_{x \in \Phi} F(x, \Phi \ominus x)} = \int_{\X} \esp{F(x, \Phi)} \sigma(\d x).
	\end{equation*}
\end{theorem}
The identity given below provides an analogous integral representation to the Campbell-Mecke formula for a PPP, 
specifically adapted to the finite-sample nature of the Binomial point process, see~\cite{DecreusefondSchulteThale2016}.
\begin{theorem}[Campbell-Mecke formula for a BPP]
	\label{th:mecke_formula_BPP}
	Let $\Phi_N$ be a BPP on $\X$ with parameter $N$ and supported by $\mu$. 
	Then, for any measurable function $F \colon \X \times \NC_{\X}^{(N)} \to \R_{+}$,
	\begin{equation*}
		\esp{\sum_{x \in \Phi_{N}} F(x, \Phi_N)} = N \int_{\X} \esp{F(x, \Phi_{N-1} \oplus x)} \mu(\d x),
	\end{equation*}
	where $\Phi_{N-1}$ is a BPP with $N-1$ points drawn i.i.d. from $\mu$.
\end{theorem}

If $\Phi$ is a PPP with intensity measure $\sigma$ and $p \in [0,1]$, then the point process denoted $p \circ \Phi$ and built by keeping with probability $p$ 
and deleting with probability $1-p$ each point $x$ of $\Phi$ independently is also a PPP with intensity measure $p \sigma$, called the $p$-thinning of $\Phi$.
In particular, a PPP $\Phi$ verifies the following invariance property:
\begin{theorem}[Invariance property]
	\label{th:invariance_property}
	For any $t \in[0,1]$,
	\begin{equation*}
		t \circ \Phi^{(1)} \oplus (1-t) \circ \Phi^{(2)} \eqdist \Phi,
	\end{equation*}
	where $\Phi^{(1)}$ and $\Phi^{(2)}$ are independent copies of $\Phi$, and $\eqdist$ denotes the equality between probability distributions.
\end{theorem}

\noindent
Using this invariance property, we can define the semigroup and the infinitesimal generator associated with a PPP through the Glauber dynamics, 
which describe a birth-death type Markov process on the configuration space.
The Glauber process evolves by randomly adding and removing points according to certain stochastic rules designed to preserve a given invariant measure.
In the case of the homogeneous PPP of intensity $\lambda$, the dynamics are particularly simple: 
each point in the configuration is removed independently at rate $1$ (death mechanism), 
and new points are added according to a homogeneous PPP of intensity $\lambda$ (birth mechanism).
This gives rise to a Markov semigroup $(\pGlaub_t)_{t \ge 0}$ acting on functions $F \in \LL^1(\NC_{\X} \to \R, \, \pi^{\sigma}) \cap \cyl(\X)$ by
\begin{equation*}
	\pGlaub_t F(\omega) = \esp{F(\G_t) \, | \, \G_0 = \omega},
\end{equation*}
where $(\G_t)_{t \ge 0}$ denotes the Glauber dynamics starting from $\omega$.
The process is ergodic with the PPP as its unique invariant measure, 
and the semigroup captures the temporal evolution of observables under this stochastic dynamics, see~\cite{EthierKurtz1986, Decreusefond2022}.

\begin{definition}[Glauber semigroup]
	\label{def:glauber_semigroup_PPP}
	Let $\Phi$ be a PPP with intensity measure $\sigma$.	We set
	\begin{equation}
		\label{eq:glauber_semigroup_PPP}
		\pGlaub_t F(\omega) = \esp{F \left( e^{-t} \circ \omega \oplus (1-e^{-t}) \circ \Phi \right)}, \; t \ge 0, \; \omega \in \NC_{\X},
	\end{equation}
	for $F \in \LL^1(\NC_{\X} \to \R, \, \pi^{\sigma})$.
\end{definition}

\begin{theorem}[Semigroup of a PPP]
	\label{th:semigroup_PPP}
	The family of operators $(\pGlaub_t)_{t \ge 0}$ given in Equation \eqref{eq:glauber_semigroup_PPP} is a semigroup on the Bannach space 
	$\LL^1(\NC_{\X} \to \R, \, \pi^{\sigma})$, whose stationary measure is the distribution of $\Phi$, i.e. 
	\begin{equation*}
		\esp{\pGlaub_t F(\Phi)} = \esp{F(\Phi)}, \; \forall \, F \in \LL^1(\NC_{\X} \to \R, \, \pi^{\sigma}).
	\end{equation*}
	Moreover, $(\pGlaub_t)_{t \ge 0}$ is an ergodic semigroup, i.e.
	\begin{equation*}
		\pGlaub_t F(\omega) \underset{t \to + \infty}{\longrightarrow} \esp{F(\Phi)}, \; \forall \, F \in \LL^1(\NC_{\X} \to \R, \, \pi^{\sigma}).
	\end{equation*}
\end{theorem}

\begin{theorem}[Infinitesimal generator of a PPP]
	\label{th:infinitesimal_generator_PPP}
	Let $\left(\pGlaub_t\right)_{t \geq 0}$ be the semigroup given in Equation \eqref{eq:glauber_semigroup_PPP}, 
	associated to the PPP $\Phi$ with intensity measure $\sigma$.
	Then, its infinitesimal generator $\lGlaub$ has its domain containing the set of measurable functions 
	$F \in \cyl(\X) \cap \Dom D$,	and is equal to 
	\begin{equation*}
		\lGlaub F(\omega) = \sum_{x \in \omega} \left(F(\omega \ominus x) - F(\omega) \right) + 
		\int_{\X} \left(F(\omega \oplus x) - F(x) \right) \sigma(\d x),
	\end{equation*}
	for $\omega \in \NC_{\X}$.
\end{theorem}

\noindent
$\Dom D$ is the domain of the discrete gradient of the PPP $\Phi$, given by:
\begin{equation*}
	\Dom D \coloneq \left\{ F \colon \NC_{\X} \to \R, \; \esp{ \int_{\X} \left| F(\Phi \oplus x) - F(\Phi) \right| \sigma(\d x)} < \infty \right\}.
\end{equation*}

\subsection{Stein's method}
In this section, we consider a topological space $E$ and we first recall the definition of $1$-Lipschitz function, which is an important notion for our work.
\begin{definition}[$1$-Lipschitz function]
	\label{def:lipschitz}
	Let $\L$ be a subset of $E$ and let $\d_{\L}$ be a distance on $\L$. 
	We say that a function $h \colon \L \to \R$ is $1$-Lipschitz according to $\d_{\L}$ if for any $x_1, x_2 \in \L$,
	\begin{equation*}
		\left|h\left(x_1\right)-h\left(x_2\right)\right| \leq \d_{\L} \left(x_1, x_2\right),
	\end{equation*}
	and we denote by $\Lip_1(\L, \d_{\L})$ the set of all these maps which are measurable.
\end{definition}

Stein's method, first introduced initially by Charles Stein in a seminal paper \cite{Stein1972} in the context of normal approximation, 
is a powerful and versatile technique developed for deriving probabilistic approximations and proving limit theorems.
More precisely, let $\eta$ be a target probability measure (e.g. Gaussian, Poisson, etc.), and suppose we wish to comapre 
a given probability distribution $\mu$ to $\eta$, both being probability distances on the same topological space $E$. 
The distances adapted to Stein's method are the ones of the form:
\begin{equation*}
	\d_{\H}(\mu, \eta) \coloneq \sup_{h \in \H} \left| \espi{\mu}{h} - \espi{\eta}{h} \right|,
\end{equation*}
where $\H$ is a space of test functions from $E$ to $\R$.

\noindent
The main metrics used in Stein's method are the Total variation distance $\DTV$ and the Wasserstein distance $\DW$ defined as follow:
\begin{definition}[Total variation distance]
	\label{def:distance_TV}
	Let $\mu$ and $\eta$ be two probability measures on a measurable space $(E, \B(E))$. We define the total variation distance on $E$ by:
	\begin{equation*}
		\begin{split}
			\DTV(\mu, \eta) &\coloneq \sup_{A \in \B(E)} \left| \mu(A) - \eta(A) \right| \\
			&= \frac{1}{2} \sup_{\| f \|_{\infty} \le 1} \left| \espi{\mu}{f} - \espi{\eta}{f} \right|.
		\end{split}
	\end{equation*}
\end{definition}
If $X, Y$ are two random variables with law $\mu$ and $\eta$ respectively, then the total variation distance writes also as:
\begin{equation*}
	\begin{split}
		\DTV(X, Y) & \coloneq \sup_{A \in \B(E)} \left| \P(X \in A) - \P(Y \in A) \right| \\
		&= \frac{1}{2} \int_E \left| f_X(z) - f_Y(z) \right| \d (\P_X + \P_Y)(z),
	\end{split}
\end{equation*}
where $f_X$ and $f_Y$ are respectively the densities of $\P_X$ and $\P_Y$ with respect to $\frac{1}{2}(\P_X + \P_Y)$. \\
In particular, if $\mu, \eta \in \FNC_{\X}$, then we have
\begin{equation*}
	\DTV(\mu, \eta) = \left| \mu \ominus \eta \right| + \left| \eta \ominus \mu \right|.
\end{equation*}

\begin{definition}[Wasserstein distance]
	\label{def:distance_W}
	Let $\mu$ and $\eta$ be two probability measures on a measurable space $(E, \B(E))$, and $\d_E$ a metric on the space $E$. 
	We define the Wasserstein distance on $E$ by:
	\begin{equation*}
		\DW(\mu, \eta) \coloneq \sup_{f \in \Lip_1(E, \d_E)} \left| \espi{\mu}{f} - \espi{\eta}{f} \right|.
	\end{equation*}
\end{definition}
Similarly, if $X$, $Y$ are two random variables with law $\mu$ and $\eta$ respectively, then the Wasserstein distance writes also as:
\begin{equation*}
	\DW(X, Y) = \sup_{f \in \Lip_1(E, \d_E)} \left| \esp{f(X)} - \esp{f(Y)} \right|.
\end{equation*}

\noindent
These metrics are well-suited to different approximation problems. 
For instance, the Wasserstein distance is particularly natural when dealing with functionals of metric spaces, 
while total variation distance is more stringent and sensitive to rare events.

\noindent
In our case, the target object is a finite Poisson point process. 
So, we consider a functional operator $\lGlaub$ which, at a formal level, satisfies for a finite point process $\Phi$ on $\X$ the identity
\begin{equation*}
	\esp{\lGlaub F(\Phi)} = 0 \text{ for a class of functions } F \colon \FNC_{\X} \to \R
\end{equation*}
if and only if, $\Phi$ is a Poisson point process with finite intensity measure $\sigma$. \\
Then, we solve the Stein's equation, which means that we ne need to find, for any test function $F \colon \FNC_{\X} \to \R$, 
a function $H_F \colon \FNC_{\X} \to \R$ such that, for any $\phi \in \FNC_{\X}$, we have
\begin{equation*}
	\lGlaub H_F(\phi) = \esp{F(\zeta)} - F(\phi),
\end{equation*}
where $\zeta$ is a Poisson point process with finite intensity measure $\sigma$.

\noindent
So, we use the so-called generator approach applied in \cite{Shih2011, CoutinDecreusefond2013} and based on theory of spatial birth-and-death processes \cite{Preston1975}. 
In our case, $\lGlaub$ is built as an infinitesimal generator of the Glauber process seen in the previous section, with the
distribution of $\zeta$ as its invariant meausre. And given the Glauber semigroup $(\pGlaub_t)_{t \ge 0}$ mentioned in Definition~\ref{def:glauber_semigroup_PPP},
we can show that we have for any $\phi \in \FNC_{\X}$,
\begin{equation*}
	\lGlaub H_F(\phi) = \int_0^{+\infty} \lGlaub \pGlaub_t F(\phi) \, \d t,
\end{equation*}
which leads us to the well-known Stein-Dirichlet representation formula for PPPs:
\begin{theorem}[Stein-Dirichlet representation formula for PPPs]
	\label{th:stein_formula_ppp}
	For any functional $F \in \cyl(\X) \cap \LL^1 (\NC_{\X} \to \R; \, \pi^{\sigma})$, and any $\omega \in \NC_{\X}$,
	the Stein-Dirichlet representation formula for the PPP $\Phi$ is given by:
	\begin{equation*}
		\esp{F(\Phi)} - F(\omega) = \int_{0}^{+\infty} \lGlaub \pGlaub_t F(\omega) \, \d t.
	\end{equation*}
\end{theorem}

\section{Quantitative Limit Theorem for Cox-Poisson Point Processes}
\label{sec:cox_PPP_limit}

\begin{figure}[!ht]
	\centering
	\includegraphics[width=7cm]{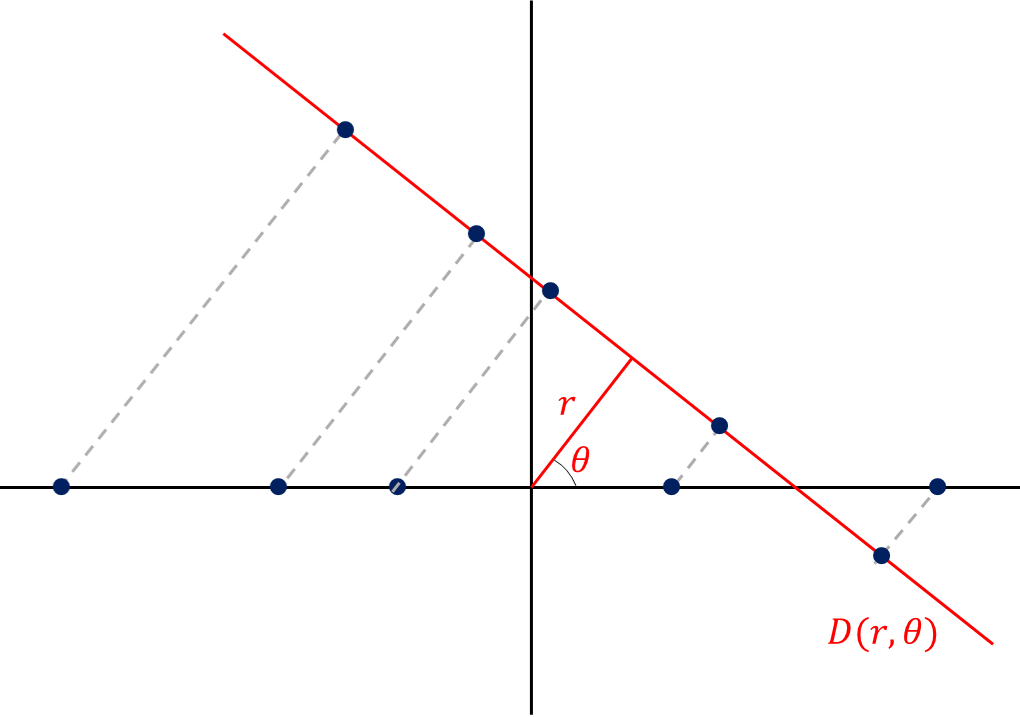}
	\caption{Cox-Poisson point process.}
	\label{fig:cox_line}
\end{figure}

The set of lines generated by a PPP on $\X \equiv \R_{+} \times[0, 2 \pi)$ is called a Poisson line process.
And a Cox-Poisson point process, as illustrated in Figure ~\ref{fig:cox_line}, is constructed by populating points on the lines of Poisson line process 
such that the locations of points on each line form an independent 1D PPP.

In other words, a Cox-Poisson point process, denoted by $Y$, can be seen as the image of a marked PPP constructed by first drawing a PPP on $\X$ 
to construct a line denoted $D(r, \theta)$ (characterized by the perpendicular distance $r$ of the line from the origin, 
and the angle $\theta$ subtended by the perpendicular dropped onto the line from the origin with respect to the positive $x$-axis in counterclockwise direction). 
And then the mark of each line is a locally finite Poisson point process on the real line $\R$.

\noindent
Our goal is to quantify the convergence of a sequence of such processes $(Y_n)_{n \in \N}$ given by the following transformation:
\begin{align*}
	T \colon \NC_{\X \times \NC_{\R}} &\longrightarrow \NC_{\R^2} \\
	\Phi_n = \displaystyle \sum_i \delta_{(r_i, \theta_i, \zeta_i^{\mu_n})} &\longmapsto 
	Y_n = \displaystyle \sum_{(r_i, \theta_i)} \left( \sum_{x \in \zeta_i^{\mu_n}} \delta_{\Pi_{(r_i, \theta_i)}(x)} \right),
\end{align*}
where $\Phi_n$ is a marked PPP on $\X \times \NC_{\R}$ with intensity 
$\sigma_n \coloneq \lambda_n \Leb^+ \otimes \frac{1}{2\pi} \carp {[0, 2\pi)} \otimes \pi^{\mu_n}$,
($\mu_n$ is the intensity of the PPP $\zeta^{\mu_n}$ on $\R$), 
and $\Pi_{(r_i, \theta_i)}(x)$ is the projection of the point $x$ from the PPP $\zeta^{\mu_n}_i$ onto the line $D(r_i, \theta_i)$.

\noindent
This analysis is rooted in the modeling principles laid out in \cite{DhillonChetlur2020}, 
and adapts technical tools from \cite{Shih2011, CoutinDecreusefond2013} to our setting in order to address the mathematical challenges 
posed by the random line structure and marked point process representation.
Following the generator approach developed in these works, 
we leverage the Stein-Dirichlet representation formula (Theorem~\ref{th:stein_formula_ppp}) as a key tool 
to bound the Wasserstein distance between the law of $Y_n$ and the target Poisson distribution.

\noindent
The following Lemma is important to prove our main result.
\begin{lemma}
	\label{lem:coarea_formula}
	Let $f \in \LL^1 \left( \R^2 \to \R_+, \; \Leb^2 \right)$. We denote by $D_{r, \theta}$ the line of $\R^2$ characterized by $(r, \theta)$.
	Then, for a fixed value of $\theta \in [0, 2\pi)$, we have:
	\begin{equation}
		\label{eq:coarea_formula_particular}
		\int_A f(x) \d x = \int_{\R_+} \left( \int_{A \cap D_{r, \theta}} f(s) \d H^1(s) \right) \d r,
		\end{equation}
	where $A$ is a subset of $\R^2$ and $H^1$ is the one-dimensional Hausdorff measure.
\end{lemma}
\noindent
Here is the main theorem of this section, stated as follows:
\begin{theorem}[Convergence of the Cox-Poisson Point Process]
	\label{th:convergence_cox_process}
	We suppose that $\lambda_n \longrightarrow + \infty \text{ and that } \mu_n = \dfrac{c}{\lambda_n}$, where $c \in \mathbb{R}^*_+$.
	We consider a compact subset $K \subset \R^2$.
	Then, we have:
	\begin{equation}
		\label{eq:convergence_rate_cox}
		\sup_{F \in \Lip_1(\NC_{K}, \DTV)} \left| \esp{F(Y_n)} - \esp{F(N)} \right| \leq 
		\frac{1}{\lambda_n} \left( c^2 \int_{\X} \left(H^1(K_{r, \theta})\right)^2 \d r \frac{\d \theta}{\pi} \right),
	\end{equation}
	where $N$ is a PPP on $\R^2$ with intensity $c \, \Leb^2$, and $K_{r, \theta} = K \cap D_{r, \theta}$.
\end{theorem}

Another way to present the previous result is to use the following inequality
\begin{equation*}
	\DW(Y_n \cap K, N \cap K) \le \frac{1}{\lambda_n} \left( c^2 \int_{\X} \left(H^1(K_{r, \theta})\right)^2 \d r \frac{\d \theta}{\pi} \right).
\end{equation*}

\section{Quantitative Limit Theorem for Cox-Binomial Point Processes}
\label{sec:cox_BPP_limit}

The BPP offers a simple yet powerful way to describe a finite set of independent points sampled according to a reference measure. 
One particularly interesting application arises in the modeling of satellites orbiting the Earth, 
where their spatial distribution can be seen as a transformation of a marked BPP.

In this setting, each satellite is initially associated with a point sampled uniformly on the unit sphere, 
representing a random position on the Earth's surface. Given such a point, we define its associated orbit as the great circle orthogonal to its position. 
This reflects the natural fact that satellites often move along nearly circular trajectories centered around the Earth. 
Furthermore, in order to capture additional randomness in their motion, we associate a mark to each point $x$.
This mark is modeled as a PPP defined on the unit circle, whose points are then projected onto the orbit of $x$ through the corresponding rotation.

Let the state space be the unit sphere $\S^2$ endowed with a reference probability measure $\nu$ (i.e. the normalized surface measure). 
We first sample a Binomial configuration of “orbits” $\{x_1,\dots,x_n\}\subset \S^2$ with i.i.d. locations from $\nu$. 
To each base point $x\in \S^2$ we associate the great circle $\Gamma_x \subset \S^2$ orthogonal to $x$ and an independent mark $\zeta_x$, 
taken to be a Poisson point process on the unit circle $\S^1$ with intensity $\mu_n$. 
We assume that $\mu_n$ represents the average number of satellites per orbit rather than per unit length. 
This implies that $\zeta(\S^1)$ follows a Poisson distribution with parameter $\mu_n$.
Writing $R_x \colon \S^1 \to \Gamma_x$ for the rotation that maps the unit circle onto $\Gamma_x$, the transformation $T$ in question is defined as follows:
\begin{align*}
	T \colon \NC_{\S^2 \times \NC_{\S^1}} &\longrightarrow \NC_{\S^2} \\
	\Phi_n = \displaystyle \sum_{i=1}^n \delta_{(x_i, \zeta_i^{\mu_n})} &\longmapsto 
	\Psi_n = \displaystyle \sum_{i=1}^n \left( \sum_{y \in \zeta_i^{\mu_n}} \delta_{R_{x_i}(y)} \right).
\end{align*}
In this setting, $\Psi_n$ is a point process on $\S^2$ representing satellite locations along great-circle orbits. 
This construction fits within the standard framework of marking and mapping theorems for point processes and their transforms; 
see, for instance, \cite{DaleyVereJones2008, MollerWaagepetersen2003, LastPenrose2017} for background, 
and \cite{MardiaJupp2000} for foundational material on statistics on spheres.

\noindent
Here is the main theorem of this section:
\begin{theorem}[Poisson approximation for the satellite model]
\label{th:convergence_satellites_BPP}
Let the marking intensity satisfy $\mu_n = \frac{c}{n}$ for some $c \in \R_{+}^{*}$, and $N$ be a PPP on $\S^2$ with intensity $c\,\nu$. 
Then, we have:
\begin{equation}
	\label{eq:convergence_satellites}
	\sup_{F \in \Lip_1(\NC_{\S^2}, \DTV)} \left| \esp{F(\Psi_n)} - \esp{F(N)} \right| \leq \frac{2c^2}{n} \cdotp
\end{equation}
\end{theorem}

\noindent
Intuitively, as the number of sampled orbits increases while individual orbits become sparsely populated, 
the superposition of independent Poisson-distributed marks leads to a spatial configuration that asymptotically exhibits the properties of a Poisson process on $\S^2$.
The intensity measure $c \, \nu$ reflects the uniform distribution of initial positions weighted by the limiting expected number of satellites per orbit.
The proof follows the Stein-generator approach adapted here to the Binomial setting. \\
In particular, the Wasserstein-$1$ distance (with respect to the underlying configuration metric of the Lipschitz class) 
between the laws of $\Psi_n$ and $N$ is at most $\frac{2c}{n}$, and hence $\Psi_n \cvdist N$ as $n\to\infty$.

\section{Conclusion}
In this paper, we have successfully derived quantitative limit theorems for two distinct classes of Cox point processes using Stein's method. 
By leveraging the generator approach and the Stein-Dirichlet representation formula for the Poisson point process, we established explicit, 
non-asymptotic bounds on the distance to a limiting PPP.

Our first result provided a convergence rate for a Cox-Poisson process on $\R^2$, constructed from a Poisson line process. 
The derived error bound highlights a clear dependence on the line intensity and the geometry of the observation window, 
offering a precise measure of approximation quality. 
Our second result focused on a Cox-Binomial model for satellite constellations on the sphere $\S^2$, for which we established a sharp convergence rate of order $O(1/n)$. 
Together, these findings demonstrate the robustness of the Stein-generator framework in handling complex dependencies 
arising from transformations of both Poisson and Binomial point processes.


\appendix
\section{Proofs}
\label{sec:proofs}

In this appendix, we provide the detailed proofs for Lemma~~\ref{lem:coarea_formula}, 
Theorems ~\ref{th:convergence_cox_process} and ~\ref{th:convergence_satellites_BPP}
\subsection{Proof of Lemma ~\protect\ref{lem:coarea_formula}}
We first recall a general reuslt known by the coarea formula given in the following Theorem (see \cite{Kuttler1998}).
\begin{theorem}[Coarea formula]
  Let $u$ be a real-valued Lipschitz function on $\R^n$, and $g \in \LL^1(\R^n \to \R_+, \, \Leb^n)$.
	The coarea formula is given by:
	\begin{equation}
    \label{eq:coarea_formula}
		\int_{\R^n} g(x) \|\nabla u(x)\| \d x = \int_{u(\R^n)} \left(\int_{u^{-1}(\{t\})} g(x) \d H^{n-1}(x) \right) \d t \, ,
	\end{equation}
	where $H^{n-1}$ is the $(n-1)$-dimensional Hausdorff measure and $\|\nabla u\|$ is the Euclidean norm of the gradient of $u$.
\end{theorem}
So, for fixed value of $r$ and $\theta$, the line $D_{r, \theta}$ is given by:
\begin{equation*}
  D_{r, \theta} = \{(x,y); \; x \cos \theta + y \sin \theta = r\}.
\end{equation*}
This motivates us to work in the particular case $n = 2$, and consider the function $u_\theta$ defined by:
\begin{equation*}
  \begin{aligned}
    u_{\theta} \colon \R^2 &\to \R_+ \\
    (x,y) &\mapsto x \cos \theta + y \sin \theta.
  \end{aligned}
\end{equation*}
This implies that $D_{r, \theta} = u_{\theta}^{-1}(\{r\})$. And since 
\begin{equation*}
  |\nabla u_{\theta}(x,y)| = \sqrt{\cos^2 \theta + \sin^2 \theta} = 1,
\end{equation*} 
then we replace in Equation ~\ref{eq:coarea_formula}, by taking the function $g \coloneq f \, \carp{A}$ ($A$ is a subset of $\R^2$), to obtain:
\begin{equation*}
		\int_{\R^2} f(x) \carp{A}(x) 1 \d x = \int_{\R_+} \left(\int_{u_{\theta}^{-1}(\{r\})} f(s) \carp{A}(s) \d H^{1}(s) \right) \d r \, ,
	\end{equation*}
which is the claimed result. \qed

\subsection{Proof of Theorem ~\protect\ref{th:convergence_cox_process}}
Fix $F\in \Lip_1(\NC_{K},\DTV)$. Then $F(\omega)=F(\omega|_{K})$ for all $\omega\in\NC_{\R^2}$.
The semigroup of the Poisson point process $N$ on $\R^2$ is given, using the Theorem ~\ref{th:semigroup_PPP}, 
for any $t \ge 0$ and $\omega \in \NC_{\R^2}$, by:
\begin{equation*}
  \pGlaub_tF(\omega) = \esp{F\left(e^{-t} \circ \omega \oplus (1-e^{-t}) \circ N' \right)},
\end{equation*}
and the infinitesimal generater of $N$ associated to $(\pGlaub_t)_{t \ge 0}$ is:
\begin{equation*}
  \lGlaub F(\omega) = \sum_{x \in \omega} \left(F(\omega \ominus x) - F(\omega) \right) + 
  c \int_{\R^2} \left(F(\omega \oplus z) - F(\omega) \right) \d z.
\end{equation*}
The Stein-Dirichlet representation formula for the Poisson point process $N$ given in Theorem ~\ref{th:stein_formula_ppp} implies that:
\begin{equation}
  \label{eq:SD}
  \esp{F(N)} - \esp{F(Y_n)} = \esp{\int_0^{+\infty} \lGlaub \pGlaub_t F(Y_n) \d t}.
\end{equation}
So, our idea for computing the convergence rate result in this theorem is to bound the expression $\left| \esp{\lGlaub F (Y_n)} \right|$,
then replacing the functional $F$ by $\pGlaub_t F$, and integrate in $t$.

\vspace{0.5em}
\step{1} \emph{Decomposition of $\esp{\lGlaub F(Y_n)}$.} \\
From the definition of $\lGlaub$ and the reduction to $K$,
\begin{multline*}
  \esp{\lGlaub F(Y_n)} = \esp{\sum_{z \in Y_n} F(Y_n \ominus z) - F(Y_n)} \\ 
    + c \, \esp{\int_K \left(F(Y_n \oplus z) - F(Y_n)\right) \d z}.
\end{multline*}
We introduce the add-subtract device along the random directions $(r,\theta) \in \Phi_n^0$ and the projected Poisson points 
$\{x_{r,\theta}\}$ onto the line $D_{r, \theta}$ (with $x_{r,\theta}=\Pi_{(r,\theta)}(x)$) to write:
\begin{multline*}
  \esp{\sum_{z \in Y_n} F(Y_n \ominus z) - F(Y_n)} = \esp{\sum_{(r,\theta) \in \Phi_n^0} \sum_{z \in K_{r, \theta}} F(Y_n \ominus z) - F(Y_n)} \\
  - \esp{\sum_{(r,\theta) \in \Phi_n^0} F\left(Y_n \ominus \sum_{x \in \zeta^{\mu_n}} \delta_{x_{r, \theta}}\right) - F(Y_n)} \\
  + \esp{\sum_{(r,\theta) \in \Phi_n^0} F\left(Y_n \ominus \sum_{x \in \zeta^{\mu_n}} \delta_{x_{r, \theta}}\right) - F(Y_n)}.
\end{multline*}

\begin{figure}[!ht]
	\centering
	\includegraphics[width=7cm]{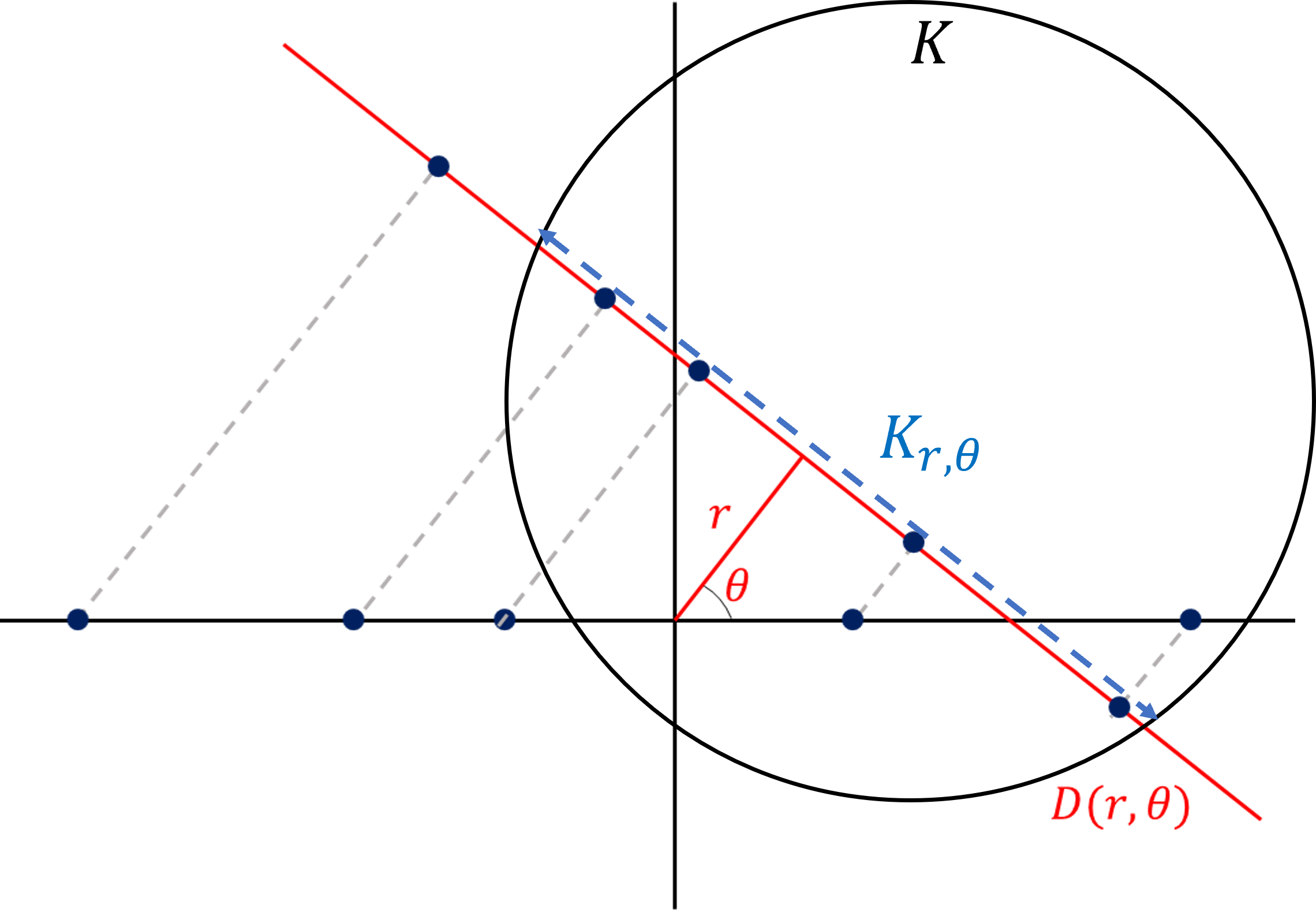}
	\caption{Reduction to $K$ of the Cox-Poisson point process.}
	\label{fig:stein_cox_line}
\end{figure}

\vspace{0.5em}
\step{2} \emph{First Campbell-Mecke step.} \\
By the Campbell-Mecke formula (Theorem~\ref{th:mecke_formula_PPP}),
\begin{multline}
  \label{eq:proof_convergence_rate_1}
  \esp{\sum_{(r,\theta) \in \Phi_n^0} F\left(Y_n \ominus \sum_{x \in \zeta^{\mu_n}} \delta_{x_{r, \theta}}\right) - F(Y_n)} \\
  = - \esp{\int_{\X \times \M} \left( F\left(Y_n \oplus \sum_{x \in \zeta^{\mu_n}} \delta_{x_{r, \theta}}\right) - F(Y_n) \right) \lambda_n \d r \frac{\d \theta}{2\pi} \pi^{\mu_n}(\d \zeta^{\mu_n})} \\
  = - \esp{\int \left( F\left(Y_n \oplus \sum_{x \in \zeta^{\mu_n}} \delta_{x_{r, \theta}}\right) - F(Y_n) \right) \carp{\zeta^{\mu_n}(K_{r, \theta}) = 1} \lambda_n \d r \frac{\d \theta}{2\pi} \pi^{\mu_n}(\d \zeta^{\mu_n})} \\
  - \esp{\int \left( F\left(Y_n \oplus \sum_{x \in \zeta^{\mu_n}} \delta_{x_{r, \theta}}\right) - F(Y_n) \right) \carp{\zeta^{\mu_n}(K_{r, \theta}) \ge 2} \lambda_n \d r \frac{\d \theta}{2\pi} \pi^{\mu_n}(\d \zeta^{\mu_n})},
\end{multline}
where $\zeta^{\mu_n}(K_{r, \theta})$ is the number of points of the Poisson process $\zeta^{\mu_n}$ projected onto the 
line $D_{r, \theta}$ that are contained in the subset $K$ (see Figure~\ref{fig:stein_cox_line}).
Which means that $\zeta^{\mu_n}(K_{r, \theta})$ has a Poisson distribution with parameter $\mu_n H^1(K_{r, \theta})$.
In the following, and for simplicity, we denote the random variable $\zeta^{\mu_n}(K_{r, \theta})$ by $U$.

\vspace{0.5em}
\step{3} \emph{Second Campbell-Mecke step and rearrangement.} \\
Applying Campbell-Mecke again and conditioning on $\{U = 0\}$, $\{U = 1\}$ and $\{U \ge 2\}$ gives:

\begin{multline*}
  \esp{\sum_{(r,\theta) \in \Phi_n^0} \sum_{z \in K_{r, \theta}} F(Y_n \ominus z) - F(Y_n)} \\
  \shoveright{- \esp{\sum_{(r,\theta) \in \Phi_n^0} F\left(Y_n \ominus \sum_{x \in \zeta^{\mu_n}} \delta_{x_{r, \theta}}\right) - F(Y_n)}} \\
  \shoveleft{= \mathbb{E} \left[ \sum_{(r,\theta) \in \Phi_n^0} \left[ \sum_{z \in K_{r, \theta}} (F(Y_n \ominus z) - F(Y_n)) \right. \right.} \\ 
  - \left. \left. \left( F\left(Y_n \ominus \sum_{x \in \zeta^{\mu_n}} \delta_{x_{r, \theta}}\right) - F(Y_n) \right) \right]\carp{U \ge 2} \right] \\
  \shoveleft{= \esp{\sum_{(r,\theta) \in \Phi_n^0} \sum_{z \in K_{r, \theta}} (F(Y_n \ominus z) - F(Y_n))\carp{U \ge 2}} } \\
  + \esp{\int_{\X \times \M} \left( F\left(Y_n \oplus \sum_{x \in \zeta^{\mu_n}} \delta_{x_{r, \theta}}\right) - F(Y_n) \right) \carp{U \ge 2} \lambda_n \d r \frac{\d \theta}{2\pi} \pi^{\mu_n}(\d \zeta^{\mu_n})},
\end{multline*}

And by adding the expression from Equation \eqref{eq:proof_convergence_rate_1}, we deduce:
\begin{multline*}
  \esp{\sum_{z \in Y_n} F(Y_n \ominus z) - F(Y_n)} = \esp{\sum_{(r,\theta) \in \Phi_n^0} \sum_{z \in K_{r, \theta}} (F(Y_n \ominus z) - F(Y_n))\carp{U \ge 2}} \\
  - \esp{\int_{\X \times \M} \left( F\left(Y_n \oplus \sum_{x \in \zeta^{\mu_n}} \delta_{x_{r, \theta}}\right) - F(Y_n) \right) \carp{U = 1} \lambda_n \d r \frac{\d \theta}{2\pi} \pi^{\mu_n}(\d \zeta^{\mu_n})} \\
  \shoveleft{= \esp{\sum_{(r,\theta) \in \Phi_n^0} \sum_{z \in K_{r, \theta}} (F(Y_n \ominus z) - F(Y_n))\carp{U \ge 2}}} \\
  - \int_{\X \times \M} \esp{\left. F\left(Y_n \oplus \sum_{x \in \zeta^{\mu_n}} \delta_{x_{r, \theta}}\right) - F(Y_n) \, \right| \, U = 1} 
  \P(U = 1) \lambda_n \d r \frac{\d \theta}{2\pi} \pi^{\mu_n}(\d \zeta^{\mu_n}) \\
  \shoveleft{= \esp{\sum_{(r,\theta) \in \Phi_n^0} \sum_{z \in K_{r, \theta}} (F(Y_n \ominus z) - F(Y_n))\carp{U \ge 2}}} \\
  - \int_{\X} \int_{K_{r,\theta}} \esp{F(Y_n \oplus s) - F(Y_n)} \frac{\d H^1(s)}{H^1(K_{r, \theta})} \mu_n H^1(K_{r, \theta}) e^{- \mu_n H^1(K_{r, \theta})} \lambda_n \d r \frac{\d \theta}{2\pi} \\
  \shoveleft{= \esp{\sum_{(r,\theta) \in \Phi_n^0} \sum_{z \in K_{r, \theta}} (F(Y_n \ominus z) - F(Y_n))\carp{U \ge 2}}} \\
  - c \, \esp{\int_{\X} \int_{K_{r, \theta}} \left(F(Y_n \oplus s) - F(Y_n)\right)\d H^1(s) e^{- \mu_n H^1(K_{r, \theta})} \d r \frac{\d \theta}{2\pi}}.
\end{multline*}
Hence, setting
\begin{equation*}
  \begin{split}
    I_1 &\coloneq c \, \esp{\int_K \left(F(Y_n \oplus z) - F(Y_n)\right) \d z} \\
    & \qquad - c \, \esp{\int_{\X} \int_{K_{r, \theta}} \left(F(Y_n \oplus s) - F(Y_n)\right)\d H^1(s) e^{- \mu_n H^1(K_{r, \theta})} \d r \frac{\d \theta}{2\pi}},
  \end{split}
\end{equation*}
and,
\begin{equation*}
  I_2 \coloneq \esp{\sum_{(r,\theta) \in \Phi_n^0} \sum_{z \in K_{r, \theta}} (F(Y_n \ominus z) - F(Y_n))\carp{\zeta^{\mu_n}(K_{r, \theta}) \ge 2}},
\end{equation*}
we obtain $\esp{\lGlaub F(Y_n)} = I_1 + I_2$ and therefore
\begin{equation*}
  \left| \esp{\lGlaub F(Y_n)} \right| \le \left| I_1 \right| + \left| I_2 \right|.
\end{equation*}

\vspace{0.5em}
\step{4} \emph{Bounds using $\Lip_1$, coarea, and $|1-e^{-x}|\le x$.}\\
$F \in \Lip_1(\NC_{K}, \DTV)$ yields $\left| F(Y_n \oplus z) - F(Y_n) \right| \le 1$.
By the coarea identity (Lemma~\ref{lem:coarea_formula}),
\begin{equation*}
  \int_K \left(F(Y_n \oplus z) - F(Y_n)\right) \d z = \int_{\X} \int_{K_{r, \theta}} \left(F(Y_n \oplus s) - F(Y_n)\right)\d H^1(s) \d r \frac{\d \theta}{2\pi} \cdotp
\end{equation*}
Hence, with $\left|1 - e^{-x} \right| \le x$ for $x \ge 0$,
\begin{equation*}
  \begin{split}
    |I_1| &\le c \, \esp{\int_{\X} \int_{K_{r, \theta}} \left| F(Y_n \oplus s) - F(Y_n) \right| \left| 1 - e^{- \mu_n H^1(K_{r, \theta})} \right| \d H^1(s) \d r \frac{\d \theta}{2 \pi}} \\
    &\le c \int_{\X} \int_{K_{r, \theta}} 1 \times \left|1 - e^{-\mu_n H^1(K_{r,\theta})} \right| \d H^1(s) \d r \frac{\d \theta}{2 \pi} \\
    &\le c \int_{\X} \left(\int_{K_{r, \theta}} \d H^1(s) \right) \mu_n H^1(K_{r,\theta})  \d r \frac{\d \theta}{2 \pi} \\
    &= c \, \mu_n \left(\int_{\X} \left(H^1(K_{r,\theta})\right)^2 \d r \frac{\d \theta}{2 \pi}\right) \\
    &= \frac{1}{\lambda_n} \left( c^2 \int_{\X} \left(H^1(K_{r,\theta})\right)^2 \d r \frac{\d \theta}{2 \pi}\right).
  \end{split}
\end{equation*}
For $I_2$, Campbell-Mecke formula gives:
\begin{equation*}
  \begin{split}
    |I_2| &\le \esp{\sum_{(r,\theta) \in \Phi_n^0} \sum_{z \in K_{r, \theta}} \carp{\zeta^{\mu_n}(K_{r, \theta}) \ge 2}} \\
    &= \int_{\X} \esp{ \zeta^{\mu_n}(K_{r, \theta}) \carp{\zeta^{\mu_n}(K_{r, \theta}) \ge 2}} \lambda_n \d r \frac{\d \theta}{2 \pi} \\
    &= \int_{\X} \mu_n H^1(K_{r, \theta})\left(1 - e^{- \mu_n H^1(K_{r, \theta})}\right) \lambda_n \d r \frac{\d \theta}{2 \pi} \\
    &\le c \int_{\X} H^1(K_{r, \theta}) \left(\mu_n H^1(K_{r, \theta}) \right) \d r \frac{\d \theta}{2 \pi} \\
    &= \frac{1}{\lambda_n} \left( c^2 \int_{\X} \left(H^1(K_{r,\theta})\right)^2 \d r \frac{\d \theta}{2 \pi}\right).
  \end{split}
\end{equation*}
Consequently,
\begin{equation}
  \label{eq:keyBound}
  \left| \esp{\lGlaub F(Y_n)} \right| \le \frac{1}{\lambda_n} \left( c^2 \int_{\X} \left(H^1(K_{r, \theta})\right)^2 \d r \frac{\d \theta}{\pi}\right).
\end{equation}

\vspace{0.5em}
\step{5} \emph{Insertion of the semigroup and integration in $t$.} \\
For $t \ge 0$,
\begin{multline*}
  \left| \pGlaub_t F(Y_n \oplus z) - \pGlaub_t F(Y_n) \right| = \left| \esp{F\left(e^{-t} \circ (Y_n \oplus z) \oplus (1-e^{-t}) \circ N \right)} \right. \\
  \shoveright{\left. - \esp{F\left(e^{-t} \circ Y_n \oplus (1-e^{-t}) \circ N \right)} \right|} \\
  \le \esp{\left| F\left(e^{-t} \circ (Y_n \oplus z) \oplus (1-e^{-t}) \circ N \right) - F\left(e^{-t} \circ Y_n \oplus (1-e^{-t}) \circ N \right) \right|} \\
  \shoveleft{= \mathbb{E} \left[ \left| F\left(e^{-t} \circ Y_n \oplus (1-e^{-t}) \circ N \oplus e^{-t} \circ z \right) \right. \right.} \\
  \shoveright{\left. \left. - F \left(e^{-t} \circ Y_n \oplus (1-e^{-t}) \circ N \right) \right| \right]} \\
  \shoveleft{= \mathbb{E} \left[ \left| F\left( \left( e^{-t} \circ Y_n \oplus (1-e^{-t}) \circ N \right) \oplus z \right) \right. \right.} \\
  \shoveright{\left. \left. - F \left(e^{-t} \circ Y_n \oplus (1-e^{-t}) \circ N \right) \right| \carp{e^{-t} \circ z = z} \right]} \\
  \le \esp{1 \times \carp{e^{-t} \circ z = z} } = \P(e^{-t} \circ z = z) = e^{-t}.
\end{multline*}
So, the bound \eqref{eq:keyBound} gains a factor $e^{-t}$ when $F$ is replaced by $\pGlaub_t F$. 
Applying \eqref{eq:SD} and integrating $e^{-t}$ over $t\in[0,\infty)$,
\begin{equation*}
  \begin{split}
    \left| \esp{F(N)} - \esp{F(Y_n)} \right| &= \left| \esp{\int_0^{+\infty} \lGlaub \pGlaub_t F(Y_n) \d t} \right| \\
    &\le \int_{0}^{+\infty} \left( e^{-t} \frac{c^2}{\lambda_n} \left(\int_{\X} \left(H^1(K_{r, \theta})\right)^2 \d r \frac{\d \theta}{\pi} \right) \right) \d t \\
    &= \left( \int_0^{+\infty} e^{-t} \right) \frac{c^2}{\lambda_n} \left(\int_{\X} \left(H^1(K_{r, \theta})\right)^2 \d r \frac{\d \theta}{\pi} \right) \\
    &= \frac{c^2}{\lambda_n} \left(\int_{\X} \left(H^1(K_{r, \theta})\right)^2 \d r \frac{\d \theta}{\pi} \right).
  \end{split}
\end{equation*}
This completes the proof. \qed

\subsection{Proof of Theorem ~\protect\ref{th:convergence_satellites_BPP}}
Let $N$ be a PPP on $\S^2$ with intensity $c\,\nu$. Its generator $\lGlaub$ is given by:
\begin{equation*}
  \lGlaub F(w) = \sum_{y \in \omega} \left(F(\omega \ominus y) - F(\omega) \right) + 
  c \int_{\S^2} \left(F(\omega \oplus y) - F(z) \right) \nu (\d y),
\end{equation*}
for $\omega \in \S^2$.
By the Stein-Dirichlet representation formula (Theorem~\ref{th:stein_formula_ppp}),
\begin{equation}
  \label{eq:stein-s2}
  \esp{F(N)} - \esp{F(\Psi_n)} = \esp{\int_0^{+\infty} \lGlaub \pGlaub_t F(\Psi_n) \d t}.
\end{equation}
Since $\pGlaub_t$, being the semigroup operator of the Poisson point process $N$, is a contraction for the add-remove metric, 
$\pGlaub_t F$ is $1$-Lipschitz whenever $F$ is. 
Hence the bounds derived below for $\lGlaub F(\Psi_n)$ apply to $\lGlaub \pGlaub_t F(\Psi_n)$ uniformly in $t$. 
It therefore suffices to bound $\left| \esp{\lGlaub F (\Psi_n)} \right|$, using the same scheme of proof of Theorem~\ref{th:convergence_cox_process}.

We write the marked Binomial point process on $\S^2\times\NC_{\S^1}$ as:
\begin{equation*}
  \Phi_n = \sum_{i=1}^n \delta_{(x_i,\zeta_i)},
\end{equation*}
with $x_i \sim \nu$ i.i.d and $\zeta_i$ is a Poisson point process with intensity $\mu_n$ on $\S^1$. 
Set $\Phi_n^0 = \displaystyle \sum_{i=1}^n \delta_{x_i}$ and 
\begin{equation*}
  \Psi_n=\sum_{i=1}^n \sum_{y \in \zeta_i}\delta_{R_{x_i}(y)}.
\end{equation*}


Expanding $\lGlaub$ at $\Psi_n$ gives:
\begin{multline}
  \label{eq:proof_convergence_satellite_1}
  \esp{\lGlaub F(\Psi_n)} = \esp{\sum_{y \in \Psi_n} F(\Psi_n \ominus y) - F(\Psi_n)} \\
  \shoveright{+ c \, \esp{\int_{\S^2} \left(F(\Psi_n \oplus y) - F(\Psi_n)\right) \nu(\d y)}} \\
  \shoveleft{= \esp{\sum_{x \in \Phi_n^0} \sum_{y \in \Gamma_x} F(\Psi_n \ominus y) - F(\Psi_n)}} \\
  - \esp{\sum_{x \in \Phi_n^0} F\left(\Psi_n \ominus \sum_{y \in \zeta^{\mu_n}} \delta_{R_x(y)} \right) - F(\Psi_n)} \\
  + \esp{\sum_{x \in \Phi_n^0} F\left(\Psi_n \ominus \sum_{y \in \zeta^{\mu_n}} \delta_{R_x(y)} \right) - F(\Psi_n)} \\
  +  c \, \esp{\int_{\S^2} \left(F(\Psi_n \oplus y) - F(\Psi_n)\right) \nu(\d y)},
\end{multline}
Let $\tilde F \colon \NC_{\S^2\times \NC_{\S^1}} \to \R$ be an intermediate functional such that
\begin{equation*}
  F(\omega) = \esp{\tilde{F}(\phi_n) \, | \, T(\phi_n) = \omega}.
\end{equation*}
It follows that
\begin{multline*}
  \esp{\sum_{x \in \Phi_n^0} F\left(\Psi_n \ominus \sum_{y \in \zeta^{\mu_n}} \delta_{R_x(y)} \right) - F(\Psi_n)} \\
  = \esp{\sum_{(x, \zeta^{\mu_n}) \in \Phi_n} \tilde{F}\left(\Phi_n \ominus \delta_{(x, \zeta^{\mu_n})} \right) - \tilde{F}(\Phi_n)}.
\end{multline*}
For simplicity, we denote the random variable $\zeta^{\mu_n}(\S^1)$ by $U$.
By the Mecke formula for Binomial processes (Theorem~\ref{th:mecke_formula_BPP}),
\begin{multline}
  \label{eq:proof_convergence_satellite_2}
  \esp{\sum_{x \in \Phi_n^0} F\left(\Psi_n \ominus \sum_{y \in \zeta^{\mu_n}} \delta_{R_x(y)} \right) - F(\Psi_n)} \\
  \shoveright{= \esp{\sum_{(x, \zeta^{\mu_n}) \in \Phi_n} \tilde{F}\left(\Phi_n \ominus \delta_{(x, \zeta^{\mu_n})} \right) - F(\Phi_n)}} \\
  \shoveleft{= n \int_{\S^2 \times \NC_{\S^1}} \esp{\tilde{F}(\Phi_{n-1}) - \tilde{F}\left(\Phi_{n-1} \oplus \delta_{(x, \zeta^{\mu_n})} \right)} \nu (\d x) \pi^{\mu_n}(\d \zeta^{\mu_n})} \\
  \shoveleft{= n \int_{\S^2 \times \NC_{\S^1}} \esp{F(\Psi_{n-1}) - F\left(\Psi_{n-1} \oplus \sum_{y \in \zeta^{\mu_n}} \delta_{R_x(y)} \right)} \nu (\d x) \pi^{\mu_n}(\d \zeta^{\mu_n})} \\
  = n \int \esp{\left( F(\Psi_{n-1}) - F\left(\Psi_{n-1} \oplus \sum_{y \in \zeta^{\mu_n}} \delta_{R_x(y)} \right) \right)\carp{U = 1}} \nu (\d x) \pi^{\mu_n}(\d \zeta^{\mu_n}) \\
  + n \int \esp{\left( F(\Psi_{n-1}) - F\left(\Psi_{n-1} \oplus \sum_{y \in \zeta^{\mu_n}} \delta_{R_x(y)} \right) \right)\carp{U \ge 2}} \nu (\d x) \pi^{\mu_n}(\d \zeta^{\mu_n}).
\end{multline}
From Equations~\eqref{eq:proof_convergence_satellite_1} and~\eqref{eq:proof_convergence_satellite_2}, we obtain:
\begin{multline*}
  \esp{\sum_{x \in \Phi_n^0} \sum_{y \in \Gamma_x} F(\Psi_n \ominus y) - F(\Psi_n)} \\
  \shoveright{- \esp{\sum_{x \in \Phi_n^0} F\left(\Psi_n \ominus \sum_{y \in \zeta^{\mu_n}} \delta_{R_x(y)} \right) - F(\Psi_n)}} \\
  \shoveleft{= \mathbb{E} \left[ \sum_{x \in \Phi_n^0} \left[ \sum_{y \in \Gamma_x}( F(\Psi_n \ominus y) - F(\Psi_n)) \right. \right.} \\ 
  - \left. \left. \left( F\left(\Psi_n \ominus \sum_{y \in \zeta^{\mu_n}} \delta_{R_x(y)} \right) - F(\Psi_n) \right) \right] \right] \\
  \shoveleft{= \mathbb{E} \left[ \sum_{x \in \Phi_n^0} \left[ \sum_{y \in \Gamma_x}( F(\Psi_n \ominus y) - F(\Psi_n)) \right. \right.} \\ 
  - \left. \left. \left( F\left(\Psi_n \ominus \sum_{y \in \zeta^{\mu_n}} \delta_{R_x(y)} \right) - F(\Psi_n) \right) \right]\carp{U \ge 2} \right] \\
  \shoveleft{= \esp{\sum_{x \in \Phi_n^0} \sum_{y \in \Gamma_x} (F(\Psi_n \ominus y) - F(\Psi_n)) \carp{U \ge 2} } } \\
  - \esp{\sum_{x \in \Phi_n^0} \left( F\left(\Psi_n \ominus \sum_{y \in \zeta^{\mu_n}} \delta_{R_x(y)} \right) - F(\Psi_n) \right) \carp{U \ge 2} }.
\end{multline*}
A second use of the Mecke formula and the use of the intermediate functional $\tilde{F}$ yield:
\begin{multline*}
  \esp{\sum_{x \in \Phi_n^0} \sum_{y \in \Gamma_x} F(\Psi_n \ominus y) - F(\Psi_n)} \\
  \shoveright{- \esp{\sum_{x \in \Phi_n^0} F\left(\Psi_n \ominus \sum_{y \in \zeta^{\mu_n}} \delta_{R_x(y)} \right) - F(\Psi_n)}} \\
  \shoveleft{= \esp{\sum_{x \in \Phi_n^0} \sum_{y \in \Gamma_x} (F(\Psi_n \ominus y) - F(\Psi_n)) \carp{U \ge 2} } } \\
  -  n \int_{\S^2 \times \NC_{\S^1}} \esp{\left( F(\Psi_{n-1}) - F\left(\Psi_{n-1} \oplus \sum_{y \in \zeta^{\mu_n}} \delta_{R_x(y)} \right) \right)\carp{U \ge 2}}
  \nu (\d x) \pi^{\mu_n}(\d \zeta^{\mu_n}).
\end{multline*}

Combining these results, we get:
\begin{multline*}
  \esp{\sum_{y \in \Psi_n} F(\Psi_n \ominus y) - F(\Psi_n)} \\
  \shoveleft{= \esp{\sum_{x \in \Phi_n^0} \sum_{y \in \Gamma_x} (F(\Psi_n \ominus y) - F(\Psi_n)) \carp{U \ge 2}}} \\
  + n \int_{\S^2 \times \NC_{\S^1}} \esp{\left( F(\Psi_{n-1}) - F\left(\Psi_{n-1} \oplus \sum_{y \in \zeta^{\mu_n}} \delta_{R_x(y)} \right) \right)\carp{U = 1}}
  \nu (\d x) \pi^{\mu_n}(\d \zeta^{\mu_n}) \\
  \shoveleft{= \esp{\sum_{x \in \Phi_n^0} \sum_{y \in \Gamma_x} (F(\Psi_n \ominus y) - F(\Psi_n)) \carp{U \ge 2}} } \\
  + n \int_{\S^2 \times \NC_{\S^1}} \esp{\left. \left( F(\Psi_{n-1}) - F\left(\Psi_{n-1} \oplus \sum_{y \in \zeta^{\mu_n}} \delta_{R_x(y)} \right) \right) \right| U = 1} \\
  \shoveright{\P(U = 1) \nu (\d x) \pi^{\mu_n}(\d \zeta^{\mu_n}) } \\
  \shoveleft{= \esp{\sum_{x \in \Phi_n^0} \sum_{y \in \Gamma_x} (F(\Psi_n \ominus y) - F(\Psi_n)) \carp{U \ge 2}} } \\
  + n \int_{\S^2} \esp{F(\Psi_{n-1}) - F(\Psi_{n-1} \oplus y)} \mu_n e^{-\mu_n} \nu(\d y) \\
  \shoveleft{= \esp{\sum_{x \in \Phi_n^0} \sum_{y \in \Gamma_x} (F(\Psi_n \ominus y) - F(\Psi_n)) \carp{U \ge 2}} } \\
  - c \int_{\S^2} \esp{F(\Psi_{n-1} \oplus y) - F(\Psi_{n-1})} e^{-\mu_n} \nu(\d y).
\end{multline*}
Therefore,
\begin{equation*}
  \esp{\lGlaub F(\Psi_n)} = I_1 + I_2 \implies \left| \esp{\lGlaub F(\Psi_n)} \right| \le \left| I_1 \right| + \left| I_2 \right|,
\end{equation*}
with
\begin{multline*}
  I_1 \coloneq c \int_{\S^2} \esp{\left(F(\Psi_n \oplus y) - F(\Psi_n)\right)} \nu(\d y) \\
  - c \int_{\S^2} \esp{F(\Psi_{n-1} \oplus y) - F(\Psi_{n-1})} e^{-\mu_n} \nu(\d y),
\end{multline*}
and,
\begin{equation*}
  I_2 \coloneq \esp{\sum_{x \in \Phi_n^0} \sum_{y \in \Gamma_x} (F(\Psi_n \ominus y) - F(\Psi_n)) \carp{\zeta^{\mu_n}(\S^1) \ge 2}}.
\end{equation*}

The difficulty we face now is that we do not have the same configuration $\Psi_n$ in the expressions of $I_1$. 
Instead, we have $\Psi_n$ and $\Psi_{n-1}$. 
To bridge the expectation over the law of $\Psi_n$ with an expression involving $\Psi_{n-1}$, 
we leverage the constructive definition of the Binomial process. 
Specifically, the $n$-orbit configuration $\Psi_n$ can be viewed as the superposition of an $(n-1)$-orbit configuration $\Psi_{n-1}$ 
and the satellite process generated from a single, independent $n$-th orbit.
This translates to the following expression:
\begin{equation*}
	\esp{F(\Psi_n)} = \esp{\int_{\S^2 \times \NC_{\S^1}} F\left(\Psi_{n-1} \oplus \sum_{z \in \zeta^{\mu_n}} \delta_{R_x(z)} \right) 
	\nu (\d x) \pi^{\mu_n}(\d \zeta^{\mu_n})}.
\end{equation*}
Therefore,
\begin{multline*}
  I_1 = c \int_{\S^2} \esp{\left(F(\Psi_n \oplus y) - F(\Psi_n)\right)} \nu(\d y) \\
  - c \int_{\S^2} \esp{F(\Psi_{n-1} \oplus y) - F(\Psi_{n-1})} e^{-\mu_n} \nu(\d y) \\
  \shoveleft{= c \int_{\S^2} \mathbb{E} \left[ \int_{\S^2 \times \NC_{\S^1}} \left( F\left(\Psi_{n-1} \oplus y \oplus \sum_{z \in \zeta^{\mu_n}} \delta_{R_x(z)} \right) \right. \right.} \\
  - \left. \left. F\left(\Psi_{n-1} \oplus \sum_{z \in \zeta^{\mu_n}} \delta_{R_x(z)} \right) \right) \nu (\d x) \pi^{\mu_n}(\d \zeta^{\mu_n}) \right] \nu(\d y) \\
  \shoveright{- c \int_{\S^2} \int_{\S^2 \times \NC_{\S^1}} \esp{F(\Psi_{n-1} \oplus y) - F(\Psi_{n-1})} e^{-\mu_n} \nu (\d x) \pi^{\mu_n}(\d \zeta^{\mu_n}) \nu(\d y)} \\
  \shoveleft{= c \int_{\S^2} \int_{\S^2 \times \NC_{\S^1}} \left[ e^{-\mu_n} \mathbb{E} \left[ \left( F\left(\Psi_{n-1} \oplus y \oplus \sum_{z \in \zeta^{\mu_n}} \delta_{R_x(z)} \right) \right. \right. \right.} \\
  - \left. \left. \left. F\left(\Psi_{n-1} \oplus \sum_{z \in \zeta^{\mu_n}} \delta_{R_x(z)} \right) \right) \right| U = 0 \right] \\
  + (1-e^{-\mu_n}) \mathbb{E} \left[ \left( F\left(\Psi_{n-1} \oplus y \oplus \sum_{z \in \zeta^{\mu_n}} \delta_{R_x(z)} \right) \right. \right.
  - \left. \left. \left. F\left(\Psi_{n-1} \oplus \sum_{z \in \zeta^{\mu_n}} \delta_{R_x(z)} \right) \right) \right| U \ge 1 \right] \\
  \shoveright{\left. - e^{-\mu_n} \esp{F(\Psi_{n-1} \oplus y) - F(\Psi_{n-1})} \right] \nu (\d x) \pi^{\mu_n}(\d \zeta^{\mu_n}) \nu(\d y)} \\
  \shoveleft{= c \int_{\S^2} \int_{\S^2 \times \NC_{\S^1}} (1- e^{-\mu_n}) \mathbb{E} \left[ \left( F\left(\Psi_{n-1} \oplus y \oplus \sum_{z \in \zeta^{\mu_n}} \delta_{R_x(z)} \right) \right. \right. } \\
  - \left. \left. \left. F\left(\Psi_{n-1} \oplus \sum_{z \in \zeta^{\mu_n}} \delta_{R_x(z)} \right) \right) \right| U \ge 1 \right] \nu (\d x) \pi^{\mu_n}(\d \zeta^{\mu_n}) \nu(\d y).
\end{multline*}
Now, using the fact that the functional $F$ is $1$-Lipschitz, which means in particular that we have
\begin{equation*}
  \left| F\left(\Psi_{n-1} \oplus y \oplus \sum_{z \in \zeta^{\mu_n}} \delta_{R_x(z)} \right) - F\left(\Psi_{n-1} \oplus \sum_{z \in \zeta^{\mu_n}} \delta_{R_x(z)} \right) \right| \le 1,
\end{equation*}
we deduce that
\begin{equation*}
  \begin{split}
    |I_1| &\le 1 c (1 - e^{- \mu_n}) \\
    &\le c \mu_n \\
    &= \frac{c^2}{n} \cdotp
  \end{split}
\end{equation*}
Again by $1$-Lipschitzness,
\begin{equation*}
  \begin{split}
    |I_2| &\le \esp{\sum_{x \in \Phi_n^0} \sum_{y \in \Gamma_x} \left| F(\Psi_n \ominus y) - F(\Psi_n) \right| \carp{\zeta^{\mu_n}(\S^1) \ge 2}} \\
    &\le \esp{\sum_{x \in \Phi_n^0} \zeta^{\mu_n}(\S^1) \carp{\zeta^{\mu_n}(\S^1) \ge 2}} \\
    &\le n \, \esp{\zeta^{\mu_n}(\S^1) \carp{\zeta^{\mu_n}(\S^1) \ge 2}} \\
    &\le n \, \mu_n (1-e^{-\mu_n}) \\
    &\le n \, \mu_n^2 \\
    &= \frac{c^2}{n} \cdotp
  \end{split}
\end{equation*}
So as a conclusion, we have shown
\begin{equation*}
  \left| \esp{\lGlaub F(\Psi_n)} \right| \leq |I_1| + |I_2| \le \frac{2 c^2}{n} \cdotp
\end{equation*}
By Equation~\eqref{eq:stein-s2} and the contraction property of $\pGlaub_t$, the same bound transfers to the integrand, yielding the claimed estimate in
\eqref{eq:convergence_satellites}. \qed
\printbibliography[heading=bibintoc]

\end{document}